\documentclass[11pt]{article}
\usepackage{mathrsfs}
\usepackage{nicefrac}
\usepackage[T1]{fontenc}
\usepackage{latexsym,amssymb,amsmath,amsfonts,amsthm}
\usepackage{graphicx}
\usepackage{caption}
\usepackage[utf8]{inputenc}
\usepackage{authblk}
\begin{document}

\centerline{\bf \Large Generalization of the Sherman--Morrison--Woodbury formula}
\centerline{\bf \Large involving the Schur complement}

\medskip

\centerline{Xuefeng Xu}

\centerline{\small Institute of Computational Mathematics and Scientific/Engineering Computing,}
\centerline{\small Academy of Mathematics and Systems Science, Chinese Academy of Sciences, Beijing 100190, China}
\centerline{\small E-mail address: xuxuefeng@lsec.cc.ac.cn}

\medskip

\noindent{\bf \large Abstract}

Let $X\in\mathbb{C}^{m\times m}$ and $Y\in\mathbb{C}^{n\times n}$ be nonsingular matrices, and let $N\in\mathbb{C}^{m\times n}$. Explicit expressions for the Moore--Penrose inverses of $M=XNY$ and a two-by-two block matrix, under appropriate conditions, have been established by Castro-Gonz\'{a}lez et al. [Linear Algebra Appl. 471 (2015) 353--368]. Based on these results, we derive a novel expression for the Moore--Penrose inverse of $A+UV^{\ast}$ under suitable conditions, where $A\in \mathbb{C}^{m\times n}$, $U\in \mathbb{C}^{m\times r}$, and $V\in \mathbb{C}^{n\times r}$. In particular, if both $A$ and $I+V^{\ast}A^{-1}U$ are nonsingular matrices, our expression reduces to the celebrated Sherman--Morrison--Woodbury formula. Moreover, we extend our results to the bounded linear operators case.

\medskip


\medskip


\medskip
\bigskip

\noindent{\bf \large 1. Introduction}

\medskip

Let $\mathbb{C}^{m\times n}$ be the set of all $m\times n$ matrices over complex field $\mathbb{C}$. The identity matrix of order $n$ is denoted by $I_{n}$ or $I$ when its size is clear in the context. For any $A\in \mathbb{C}^{m\times n}$, let $A^{\ast}$, $\mathcal{R}(A)$, and $\mathcal{N}(A)$ denote the conjugate transpose, the range, and the null space of $A$, respectively. The \emph{Moore--Penrose (MP) inverse} of $A\in \mathbb{C}^{m\times n}$ is denoted by $A^{\dagger}$, which is defined as the unique matrix $Z\in\mathbb{C}^{n\times m}$ satisfying the following equations:
\begin{displaymath}
(1)\ AZA=A,\quad (2)\ ZAZ=Z,\quad (3)\ (AZ)^{\ast}=AZ,\quad (4)\ (ZA)^{\ast}=ZA.
\end{displaymath}
Clearly, the MP inverse $A^{\dagger}$ coincides with the usual inverse $A^{-1}$ when $A$ is nonsingular. The symbols $E_{A}=I-AA^{\dagger}$ and $F_{A}=I-A^{\dagger}A$ denote the orthogonal projectors onto $\mathcal{N}(A^{\ast})$ and $\mathcal{N}(A)$, respectively. A matrix $Z\in\mathbb{C}^{n\times m}$ is referred to as a \emph{$\{1\}$-inverse} of $A$ if it satisfies the equality~(1); see, e.g., [1, Chapter 1, Definition 1].

Let $A\in\mathbb{C}^{n\times n}$, $U\in\mathbb{C}^{n\times r}$, and $V\in\mathbb{C}^{n\times r}$. If both $A$ and $I+V^{\ast}A^{-1}U$ are nonsingular, then $A+UV^{\ast}$ is also nonsingular and
\begin{align*}
\left(A+UV^{\ast}\right)^{-1}=A^{-1}-A^{-1}U\left(I+V^{\ast}A^{-1}U\right)^{-1}V^{\ast}A^{-1},\tag{1.1}
\end{align*}
which is the celebrated \emph{Sherman--Morrison--Woodbury (SMW) formula} (see [2--4]). Assume that $A^{-1}$ has been precomputed. If $r$ is much smaller than $n$, then $I+V^{\ast}A^{-1}U$ is much easier to invert than $A+UV^{\ast}$. Hence, the formula (1.1) provides an effective way to compute $\left(A+UV^{\ast}\right)^{-1}$. 

The SMW formula is widely used in many fields, such as statistics, networks, structural analysis, asymptotic analysis, optimization and partial differential equations; see, e.g., [5]. Here we only mention two specific applications of the SMW formula. Using the SMW formula, Malyshev and Sadkane [6] obtained a fast numerical algorithm for solving systems of linear equations with tridiagonal block Toeplitz matrices. Lai and Vemuri [7] applied the SMW formula to solve the surface smoothing problem. Using finite element methods to discretize the variational formulation of the surface smoothing problem may yield a linear system. The SMW formula can convert the problem of solving the original linear system to solving a Lyapunov matrix equation or a cascade of two Lyapunov matrix equations. The simplified problem can be solved efficiently using the ADI method and the bi-conjugate-gradient technique.

However, the SMW formula (1.1) is invalid when $A$ is not square or the assumption is not satisfied. Suppose that $A\in\mathbb{C}^{m\times n}$ and $b\in\mathbb{C}^{m}$ (here $\mathbb{C}^{m}=\mathbb{C}^{m\times 1}$). It is well-known that the minimum $\ell_{2}$-norm solution of the Least Squares Problem (LSP)
\begin{displaymath}
\min_{x\in\mathbb{C}^{n}}\|Ax-b\|_{2}
\end{displaymath}
is $x_{\ast}=A^{\dagger}b$. This solution vector can be statistically interpreted as providing an optimal estimator among all linear unbiased estimators, and it can be geometrically interpreted as providing an orthogonal projection of $b$ onto $\mathcal{R}(A)$. Because of the relationship with the LSP, the MP inverse may be the most important of all other generalized
inverses. The perturbation theory of the MP inverse is a classical topic in matrix analysis and numerical linear algebra; see, e.g., [8--11].

In general, the perturbations of a matrix can be divided into two categories: additive type and multiplicative type. In this paper, we derive an explicit expression for the MP inverse of the additive perturbations of a matrix using the results in [12, Theorems 2.2 and 3.2], which generalizes the classical SMW formula (1.1). Let $A\in \mathbb{C}^{m\times n}$, $U\in \mathbb{C}^{m\times r}$, and $V\in \mathbb{C}^{n\times r}$. If $\mathcal{R}(U)\subseteq \mathcal{R}(A)$, $\mathcal{R}(V)\subseteq \mathcal{R}(A^{\ast})$, $\mathcal{R}(U^{\ast})\subseteq \mathcal{R}(S_{A})$, and $\mathcal{R}(V^{\ast})\subseteq \mathcal{R}(S_{A}^{\ast})$, then we have
\begin{align*}
(A+UV^{\ast})^{\dagger}=\big(I+A^{\dagger}UF_{S_{A}}U^{\ast}(A^{\dagger})^{\ast}\big)^{-1}\big(A^{\dagger}-A^{\dagger}US_{A}^{\dagger}V^{\ast}A^{\dagger}\big)\big(I+(A^{\dagger})^{\ast}VE_{S_{A}}V^{\ast}A^{\dagger}\big)^{-1},\tag{1.2}
\end{align*}
where $S_{A}=I+V^{\ast}A^{\dagger}U$. Note that, if both $A$ and $I+V^{\ast}A^{-1}U$ are nonsingular, the conditions are automatically satisfied and our expression (1.2) reduces to (1.1).

The rest of this paper is organized as follows. In Section 2, we present two useful lemmas, which play an important role in our subsequent derivations. In Section 3, we give an explicit expression for $(A+UV^{\ast})^{\dagger}$, provided that $\mathcal{R}(U)\subseteq \mathcal{R}(A)$, $\mathcal{R}(V)\subseteq \mathcal{R}(A^{\ast})$, $\mathcal{R}(U^{\ast})\subseteq \mathcal{R}(S_{A})$, and $\mathcal{R}(V^{\ast})\subseteq \mathcal{R}(S_{A}^{\ast})$. In Section 4, we extend the results established in Theorems 3.2 and 3.7 below to the bounded linear operators case.

\medskip
\bigskip

\noindent{\bf \large 2. Preliminaries}

\medskip

In order to prove the expression (1.2), we need the following two lemmas (see [12, Theorems 2.2 and 3.2]), which play a key role in our subsequent derivations.

\medskip

Let $X\in \mathbb{C}^{m\times m}$ and $Y\in\mathbb{C}^{n\times n}$ be nonsingular matrices, and let $N\in \mathbb{C}^{m\times n}$. The following lemma gives an explicit expression for $(XNY)^{\dagger}$, provided that $XE_{N}=E_{N}$ and $F_{N}Y=F_{N}$.

\medskip

\noindent{\bf Lemma 2.1.} \emph{Let $N\in \mathbb{C}^{m\times n}$, $X\in \mathbb{C}^{m\times m}$, $Y\in\mathbb{C}^{n\times n}$, and $M=XNY$. If both $X$ and $Y$ are nonsingular, $XE_{N}=E_{N}$, and $F_{N}Y=F_{N}$, then}
\begin{align*}
M^{\dagger}=(I+L^{\ast})(I+LL^{\ast})^{-1}Y^{-1}N^{\dagger}X^{-1}(I+R^{\ast}R)^{-1}(I+R^{\ast}),\tag{2.1}
\end{align*}
\emph{where $R=E_{N}(I-X^{-1})$ and $L=(I-Y^{-1})F_{N}$.}

\medskip

Generalized inverses of partitioned matrices possess some important and interesting properties; see, e.g., [13--16]. Now, let $M$ be the following two-by-two block matrix
\begin{align*}
M=\begin{pmatrix}
A&C\\
B&D\\
\end{pmatrix},\tag{2.2}
\end{align*}
where $A\in\mathbb{C}^{p\times q}$, $B\in\mathbb{C}^{r\times q}$, $C\in\mathbb{C}^{p\times s}$, and $D\in\mathbb{C}^{r\times s}$ ($p$, $q$, $r$, and $s$ are all positive integers) are the corresponding submatrices of $M$. The matrix $S_{A}=D-BA^{\dagger}C$
is called the \emph{generalized Schur complement} of $A$ in $M$. If $\mathcal{R}(B^{\ast})\subseteq \mathcal{R}(A^{\ast})$ and $\mathcal{R}(C)\subseteq \mathcal{R}(A)$, then we have the following expression for $M^{\dagger}$.

\medskip

\noindent{\bf Lemma 2.2.} \emph{Let $M$ be a block matrix of the form {\rm (2.2)}. If $\mathcal{R}(B^{\ast})\subseteq \mathcal{R}(A^{\ast})$ and $\mathcal{R}(C)\subseteq \mathcal{R}(A)$, then}
\begin{align*}
M^{\dagger}=\begin{pmatrix}
\Sigma&\Sigma H^{\ast}E_{S_{A}}-\Psi KS_{A}^{\dagger}\\
F_{S_{A}}K^{\ast}\Sigma-S_{A}^{\dagger}H\Phi&S_{A}^{\dagger}-S_{A}^{\dagger}H\Phi H^{\ast}E_{S_{A}}-F_{S_{A}}K^{\ast}\Psi KS_{A}^{\dagger}+F_{S_{A}}K^{\ast}\Sigma H^{\ast}E_{S_{A}}
\end{pmatrix},\tag{2.3}
\end{align*}
\emph{where} 
\begin{displaymath}
H=BA^{\dagger},\ K=A^{\dagger}C, \ \Phi=(I+H^{\ast}E_{S_{A}}H)^{-1}, \ \Psi=(I+KF_{S_{A}}K^{\ast})^{-1}, \ \Sigma=\Psi(A^{\dagger}+KS_{A}^{\dagger}H)\Phi.
\end{displaymath}

\bigskip

\noindent{\bf \large 3. Main results}

\medskip

In order to prove our main formula in Theorem 3.2 below, we first give an important lemma.

\medskip

\noindent{\bf Lemma 3.1.} \emph{Let $A\in \mathbb{C}^{m\times n}$, $U\in \mathbb{C}^{m\times r}$, $V\in \mathbb{C}^{n\times r}$,
$X=\begin{pmatrix}
I&-U\\
0&I\\
\end{pmatrix}$, $N=\begin{pmatrix}
A&U\\
-V^{\ast}&I\\
\end{pmatrix}$, and $Y=\begin{pmatrix}
I&0\\
V^{\ast}&I
\end{pmatrix}$.
If $\mathcal{R}(U)\subseteq \mathcal{R}(A)$ and  $\mathcal{R}(V)\subseteq \mathcal{R}(A^{\ast})$, then} 
\begin{displaymath}
XE_{N}=E_{N}\Longleftrightarrow \mathcal{R}(U^{\ast})\subseteq \mathcal{R}(S_{A}),\quad F_{N}Y=F_{N}\Longleftrightarrow \mathcal{R}(V^{\ast})\subseteq \mathcal{R}(S_{A}^{\ast}),
\end{displaymath}
\emph{where $S_{A}=I+V^{\ast}A^{\dagger}U$.}

\medskip

\noindent{\bf Proof.} (i) Due to $\mathcal{R}(U)\subseteq \mathcal{R}(A)$ and  $\mathcal{R}(V)\subseteq \mathcal{R}(A^{\ast})$, it follows from Lemma 2.2 that
\begin{displaymath}
N^{\dagger}=\begin{pmatrix}
N_{1} & N_{3}\\
N_{2} & N_{4}\\
\end{pmatrix},
\end{displaymath}
where
\begin{align*}
N_{1}&=\Sigma,\tag{3.1a}\\ N_{2}&=F_{S_{A}}K^{\ast}\Sigma-S_{A}^{\dagger}H\Phi,\tag{3.1b}\\ 
N_{3}&=\Sigma H^{\ast}E_{S_{A}}-\Psi KS_{A}^{\dagger},\tag{3.1c}\\
N_{4}&=S_{A}^{\dagger}-S_{A}^{\dagger}H\Phi H^{\ast}E_{S_{A}}-F_{S_{A}}K^{\ast}\Psi KS_{A}^{\dagger}+F_{S_{A}}K^{\ast}\Sigma H^{\ast}E_{S_{A}}.\tag{3.1d}
\end{align*}
And here
\begin{align*}
S_{A}=I+V^{\ast}&A^{\dagger}U, \quad H=-V^{\ast}A^{\dagger}, \quad K=A^{\dagger}U,\tag{3.2a}\\
\Phi=(I+H^{\ast}E_{S_{A}}H)^{-1},\quad \Psi&=(I+KF_{S_{A}}K^{\ast})^{-1}, \quad \Sigma=\Psi(A^{\dagger}+KS_{A}^{\dagger}H)\Phi.\tag{3.2b}
\end{align*}
Let $E_{N}=\begin{pmatrix}
E_{1}&E_{3}\\
E_{2}&E_{4}\\
\end{pmatrix}$, $XE_{N}=\begin{pmatrix}
\widetilde{E}_{1}&\widetilde{E}_{3}\\
\widetilde{E}_{2}&\widetilde{E}_{4}\\
\end{pmatrix}$, $F_{N}=\begin{pmatrix}
F_{1}&F_{3}\\
F_{2}&F_{4}\\
\end{pmatrix}$, and $F_{N}Y=\begin{pmatrix}
\widetilde{F}_{1}&\widetilde{F}_{3}\\
\widetilde{F}_{2}&\widetilde{F}_{4}\\
\end{pmatrix}$. Direct computations yield
\begin{align*}
E_{1}&=I-AN_{1}-UN_{2}, \ 
E_{2}=V^{\ast}N_{1}-N_{2},\ 
E_{3}=-AN_{3}-UN_{4}, \ 
E_{4}=I+V^{\ast}N_{3}-N_{4},\tag{3.3a}\\
\widetilde{E}_{1}&=I-(A+UV^{\ast})N_{1},\ 
\widetilde{E}_{2}=E_{2},\ \widetilde{E}_{3}=-U-(A+UV^{\ast})N_{3}, \
\widetilde{E}_{4}=E_{4},\tag{3.3b}\\
F_{1}&=I-N_{1}A+N_{3}V^{\ast}, \
F_{2}=-N_{2}A+N_{4}V^{\ast},\ F_{3}=-N_{1}U-N_{3}, \ 
F_{4}=I-N_{2}U-N_{4},\tag{3.3c}\\
\widetilde{F}_{1}&=I-N_{1}(A+UV^{\ast}), \
\widetilde{F}_{2}=V^{\ast}-N_{2}(A+UV^{\ast}),\ 
\widetilde{F}_{3}=F_{3},\ 
\widetilde{F}_{4}=F_{4}.\tag{3.3d}
\end{align*}

(ii) We claim that the following equalities hold:
\begin{align*}
E_{1}-\widetilde{E}_{1}&=UE_{S_{A}}V^{\ast}A^{\dagger}\Phi,\qquad E_{3}-\widetilde{E}_{3}=UE_{S_{A}}(I+HH^{\ast}E_{S_{A}})^{-1},\tag{3.4a}\\
F_{1}-\widetilde{F}_{1}&=\Psi A^{\dagger}UF_{S_{A}}V^{\ast},\ \qquad F_{2}-\widetilde{F}_{2}=-(I+F_{S_{A}}K^{\ast} K)^{-1}F_{S_{A}}V^{\ast}.\tag{3.4b}
\end{align*}
In fact,
\begin{align*}
E_{1}-\widetilde{E}_{1}&=U(V^{\ast}N_{1}-N_{2})\\
&=U\big[V^{\ast}\Psi(A^{\dagger}+KS_{A}^{\dagger}H)-F_{S_{A}}K^{\ast}\Psi(A^{\dagger}+KS_{A}^{\dagger}H)+S_{A}^{\dagger}H\big]\Phi\\
&=U\big[V^{\ast}\Psi A^{\dagger}+V^{\ast}\Psi KS_{A}^{\dagger}H-F_{S_{A}}K^{\ast}\Psi A^{\dagger}-F_{S_{A}}K^{\ast}\Psi KS_{A}^{\dagger}H+S_{A}^{\dagger}H\big]\Phi\\
&=U\big[V^{\ast}\Psi A^{\dagger}+V^{\ast}\Psi KS_{A}^{\dagger}H-F_{S_{A}}K^{\ast}\Psi A^{\dagger}+(I-F_{S_{A}}K^{\ast}\Psi K)S_{A}^{\dagger}H\big]\Phi\\
&=U\big[V^{\ast}\Psi A^{\dagger}+V^{\ast}\Psi KS_{A}^{\dagger}H-F_{S_{A}}K^{\ast}\Psi A^{\dagger}+(I+F_{S_{A}}K^{\ast}K)^{-1}S_{A}^{\dagger}H\big]\Phi, \tag{3.5}
\end{align*}
where we used the fact that $I-F_{S_{A}}K^{\ast}\Psi K=(I+F_{S_{A}}K^{\ast}K)^{-1}$. By formula (1.1), we have
\begin{displaymath}
\Psi=(I+KF_{S_{A}}K^{\ast})^{-1}=I-K(I+F_{S_{A}}K^{\ast}K)^{-1}F_{S_{A}}K^{\ast},
\end{displaymath}
which yields
\begin{align*}
V^{\ast}\Psi A^{\dagger}&=V^{\ast}A^{\dagger}-V^{\ast}K(I+F_{S_{A}}K^{\ast}K)^{-1}F_{S_{A}}K^{\ast}A^{\dagger}, \tag{3.6}\\
V^{\ast}\Psi KS_{A}^{\dagger}H&=V^{\ast}KS_{A}^{\dagger}H-V^{\ast}K(I+F_{S_{A}}K^{\ast}K)^{-1}F_{S_{A}}K^{\ast}KS_{A}^{\dagger}H,
\tag{3.7}\\
F_{S_{A}}K^{\ast}\Psi A^{\dagger}&=F_{S_{A}}K^{\ast}A^{\dagger}-F_{S_{A}}K^{\ast}K(I+F_{S_{A}}K^{\ast}K)^{-1}F_{S_{A}}K^{\ast}A^{\dagger}.\tag{3.8}
\end{align*}
Note that $V^{\ast}K=S_{A}-I$. By substituting $V^{\ast}K=S_{A}-I$ into (3.6) and (3.7), we obtain
\begin{align*}
V^{\ast}\Psi A^{\dagger}&=V^{\ast}A^{\dagger}-(S_{A}-I)(I+F_{S_{A}}K^{\ast}K)^{-1}F_{S_{A}}K^{\ast}A^{\dagger},\\
V^{\ast}\Psi KS_{A}^{\dagger}H&=(S_{A}-I)S_{A}^{\dagger}H-(S_{A}-I)(I+F_{S_{A}}K^{\ast}K)^{-1}F_{S_{A}}K^{\ast}KS_{A}^{\dagger}H.
\end{align*}
Due to $S_{A}F_{S_{A}}=S_{A}(I-S_{A}^{\dagger}S_{A})=0$, it follows that 
\begin{displaymath}
-(S_{A}-I)(I+F_{S_{A}}K^{\ast}K)^{-1}F_{S_{A}}=(I+F_{S_{A}}K^{\ast}K)^{-1}F_{S_{A}}.
\end{displaymath}
Hence,
\begin{align*}
V^{\ast}\Psi A^{\dagger}&=V^{\ast}A^{\dagger}+(I+F_{S_{A}}K^{\ast}K)^{-1}F_{S_{A}}K^{\ast}A^{\dagger},\tag{3.9}\\
V^{\ast}\Psi KS_{A}^{\dagger}H&=(S_{A}-I)S_{A}^{\dagger}H+(I+F_{S_{A}}K^{\ast}K)^{-1}F_{S_{A}}K^{\ast}KS_{A}^{\dagger}H.\tag{3.10}
\end{align*}
From (3.8) and (3.9), we have
\begin{align*}
V^{\ast}\Psi A^{\dagger}-F_{S_{A}}K^{\ast}\Psi A^{\dagger}=V^{\ast}A^{\dagger}-F_{S_{A}}K^{\ast}A^{\dagger}+F_{S_{A}}K^{\ast}A^{\dagger}=V^{\ast}A^{\dagger}.\tag{3.11}
\end{align*}
From (3.10), we have
\begin{align*}
V^{\ast}\Psi KS_{A}^{\dagger}H+(I+F_{S_{A}}K^{\ast}K)^{-1}S_{A}^{\dagger}H=(S_{A}-I)S_{A}^{\dagger}H+S_{A}^{\dagger}H=S_{A}S_{A}^{\dagger}H.\tag{3.12}
\end{align*}
Inserting (3.11) and (3.12) into (3.5) gives
\begin{displaymath}
E_{1}-\widetilde{E}_{1}=U(V^{\ast}A^{\dagger}+S_{A}S_{A}^{\dagger}H)\Phi=U(V^{\ast}A^{\dagger}-S_{A}S_{A}^{\dagger}V^{\ast}A^{\dagger})\Phi=UE_{S_{A}}V^{\ast}A^{\dagger}\Phi.
\end{displaymath}

We next verify the second equality in (3.4a). On account of (1.1), (3.1c), (3.1d), (3.2a), and (3.2b), we obtain
\begin{align*}
N_{4}&=S_{A}^{\dagger}-S_{A}^{\dagger}H\Phi H^{\ast}E_{S_{A}}-F_{S_{A}}K^{\ast}\Psi KS_{A}^{\dagger}+F_{S_{A}}K^{\ast}\Sigma H^{\ast}E_{S_{A}}\\
&=S_{A}^{\dagger}(I+HH^{\ast}E_{S_{A}})^{-1}-F_{S_{A}}K^{\ast}\Psi KS_{A}^{\dagger}+F_{S_{A}}K^{\ast}\Psi(A^{\dagger}+KS_{A}^{\dagger}H)\Phi H^{\ast}E_{S_{A}}\\
&=S_{A}^{\dagger}(I+HH^{\ast}E_{S_{A}})^{-1}+F_{S_{A}}K^{\ast}\Psi A^{\dagger}\Phi H^{\ast}E_{S_{A}}-F_{S_{A}}K^{\ast}\Psi KS_{A}^{\dagger}(I+HH^{\ast}E_{S_{A}})^{-1}\\
&=(I+F_{S_{A}}K^{\ast}K)^{-1}S_{A}^{\dagger}(I+HH^{\ast}E_{S_{A}})^{-1}+F_{S_{A}}K^{\ast}\Psi A^{\dagger}\Phi H^{\ast}E_{S_{A}},\tag{3.13}\\
V^{\ast}N_{3}&=V^{\ast}\Sigma H^{\ast}E_{S_{A}}-V^{\ast}\Psi KS_{A}^{\dagger}\\
&=V^{\ast}\Psi(A^{\dagger}+KS_{A}^{\dagger}H)\Phi H^{\ast}E_{S_{A}}-V^{\ast}\Psi KS_{A}^{\dagger}\\
&=V^{\ast}\Psi A^{\dagger}\Phi H^{\ast}E_{S_{A}}-V^{\ast}\Psi KS_{A}^{\dagger}(I+HH^{\ast}E_{S_{A}})^{-1},\tag{3.14}\\
\Phi H^{\ast}E_{S_{A}}&=(I+H^{\ast}E_{S_{A}}H)^{-1}H^{\ast}E_{S_{A}}\\
&=H^{\ast}E_{S_{A}}-H^{\ast}E_{S_{A}}(I+HH^{\ast}E_{S_{A}})^{-1}HH^{\ast}E_{S_{A}}\\
&=H^{\ast}E_{S_{A}}(I+HH^{\ast}E_{S_{A}})^{-1}.\tag{3.15}
\end{align*}
In view of (3.3a), (3.3b), (3.13), (3.14), and (3.15), we derive
\begin{align*}
E_{3}-\widetilde{E}_{3}=U(I-N_{4}+V^{\ast}N_{3})=U\widehat{E}_{3}(I+HH^{\ast}E_{S_{A}})^{-1},\tag{3.16}
\end{align*}
where
\begin{displaymath}
\widehat{E}_{3}=I+HH^{\ast}E_{S_{A}}-(I+F_{S_{A}}K^{\ast}K)^{-1}S_{A}^{\dagger}+(V^{\ast}\Psi -F_{S_{A}}K^{\ast}\Psi )A^{\dagger}H^{\ast}E_{S_{A}}-V^{\ast}\Psi KS_{A}^{\dagger}.
\end{displaymath}
Using $\Psi=I-K(I+F_{S_{A}}K^{\ast}K)^{-1}F_{S_{A}}K^{\ast}$, $V^{\ast}K=S_{A}-I$, and $S_{A}F_{S_{A}}=0$, we obtain
\begin{align*}
V^{\ast}\Psi&=V^{\ast}-(S_{A}-I)(I+F_{S_{A}}K^{\ast}K)^{-1}F_{S_{A}}K^{\ast}=V^{\ast}+(I+F_{S_{A}}K^{\ast}K)^{-1}F_{S_{A}}K^{\ast},\\
F_{S_{A}}K^{\ast}\Psi&=F_{S_{A}}K^{\ast}-F_{S_{A}}K^{\ast}K(I+F_{S_{A}}K^{\ast}K)^{-1}F_{S_{A}}K^{\ast}=(I+F_{S_{A}}K^{\ast}K)^{-1}F_{S_{A}}K^{\ast},\\
V^{\ast}\Psi KS_{A}^{\dagger}&=(S_{A}-I)S_{A}^{\dagger}+(I+F_{S_{A}}K^{\ast}K)^{-1}F_{S_{A}}K^{\ast}KS_{A}^{\dagger}=S_{A}S_{A}^{\dagger}-(I+F_{S_{A}}K^{\ast}K)^{-1}S_{A}^{\dagger}.
\end{align*}
Hence,
\begin{align*}
\widehat{E}_{3}&=I+HH^{\ast}E_{S_{A}}-(I+F_{S_{A}}K^{\ast}K)^{-1}S_{A}^{\dagger}-HH^{\ast}E_{S_{A}}-S_{A}S_{A}^{\dagger}+(I+F_{S_{A}}K^{\ast}K)^{-1}S_{A}^{\dagger}\\
&=E_{S_{A}}.\tag{3.17}
\end{align*}
By plugging (3.17) back into (3.16), we get that $E_{3}-\widetilde{E}_{3}=UE_{S_{A}}(I+HH^{\ast}E_{S_{A}})^{-1}$.

Analogously, we can show that the equalities in (3.4b) also hold (note that $HU=I-S_{A}$ and $E_{S_{A}}S_{A}=0$).

(iii) If $XE_{N}=E_{N}$, then $UE_{S_{A}}=0$, or equivalently, $\mathcal{R}(U^{\ast})\subseteq \mathcal{R}(S_{A})$. Conversely, if $\mathcal{R}(U^{\ast})\subseteq \mathcal{R}(S_{A})$, then we conclude that $XE_{N}=E_{N}$ due to (3.4a). Similarly, we can verify that $F_{N}Y=F_{N}$ if and only if  $\mathcal{R}(V^{\ast})\subseteq \mathcal{R}(S_{A}^{\ast})$. This completes the proof. \qed

\medskip

Based on Lemmas 2.1 and 3.1, we can prove the following result, which provides an explicit expression for $(A+UV^{\ast})^{\dagger}$ under suitable conditions.

\medskip

\noindent{\bf Theorem 3.2.} \emph{Let $A\in \mathbb{C}^{m\times n}$, $U\in \mathbb{C}^{m\times r}$, and $V\in \mathbb{C}^{n\times r}$. If $\mathcal{R}(U)\subseteq \mathcal{R}(A)$, $\mathcal{R}(V)\subseteq \mathcal{R}(A^{\ast})$, $\mathcal{R}(U^{\ast})\subseteq \mathcal{R}(S_{A})$, and $\mathcal{R}(V^{\ast})\subseteq \mathcal{R}(S_{A}^{\ast})$, then the MP inverse of $A+UV^{\ast}$ is given by
\begin{align*}
(A+UV^{\ast})^{\dagger}=\big(I+A^{\dagger}UF_{S_{A}}U^{\ast}(A^{\dagger})^{\ast}\big)^{-1}\big(A^{\dagger}-A^{\dagger}US_{A}^{\dagger}V^{\ast}A^{\dagger}\big)\big(I+(A^{\dagger})^{\ast}VE_{S_{A}}V^{\ast}A^{\dagger}\big)^{-1},\tag{3.18}
\end{align*}
where $S_{A}=I+V^{\ast}A^{\dagger}U$.} 

\medskip

\noindent{\bf Proof.} 
(i) We first note that
\begin{align*}
\mathcal{R}(U)\subseteq \mathcal{R}(A)&\Longleftrightarrow E_{A}U=0,\qquad \mathcal{R}(V)\subseteq \mathcal{R}(A^{\ast})\Longleftrightarrow V^{\ast}F_{A}=0,\\
\mathcal{R}(U^{\ast})\subseteq \mathcal{R}(S_{A})&\Longleftrightarrow UE_{S_{A}}=0,\quad  \mathcal{R}(V^{\ast})\subseteq \mathcal{R}(S_{A}^{\ast})\Longleftrightarrow F_{S_{A}}V^{\ast}=0.
\end{align*}
Let $X$, $N$, and $Y$ be defined as in Lemma 3.1, and let $M=XNY$. The assumption of this theorem implies that $XE_{N}=E_{N}$ and $F_{N}Y=F_{N}$ due to Lemma 3.1. Therefore, we can apply formula (2.1) to compute $M^{\dagger}$. 

We now calculate some quantities involved in (2.1). Straightforward  computations yield
\begin{displaymath}
R=E_{N}(I-X^{-1})=\begin{pmatrix}
0&-E_{1}U\\
0&-E_{2}U\\
\end{pmatrix},\quad L=(I-Y^{-1})F_{N}=\begin{pmatrix}
0&0\\
V^{\ast}F_{1}&V^{\ast}F_{3}\\
\end{pmatrix}.
\end{displaymath}
Hence,
\begin{align*}
&I+L^{\ast}=\begin{pmatrix}
I&F_{1}^{\ast}V\\
0&I+F_{3}^{\ast}V\\
\end{pmatrix},\quad  I+R^{\ast}=\begin{pmatrix}
I&0\\
-U^{\ast}E_{1}^{\ast}&I-U^{\ast}E_{2}^{\ast}\\
\end{pmatrix},\\
(I&+LL^{\ast})^{-1}=\begin{pmatrix}
I&0\\
0&(I+V^{\ast}F_{1}F_{1}^{\ast}V+V^{\ast}F_{3}F_{3}^{\ast}V)^{-1}\\
\end{pmatrix}=\begin{pmatrix}
I & 0\\
0 & \widehat{F}\\
\end{pmatrix},\\
(I&+R^{\ast}R)^{-1}=\begin{pmatrix}
I&0\\
0&(I+U^{\ast}E_{1}^{\ast}E_{1}U+U^{\ast}E_{2}^{\ast}E_{2}U)^{-1}\\
\end{pmatrix}=\begin{pmatrix}
I & 0\\
0 & \widehat{E}\\
\end{pmatrix},
\end{align*}
where 
\begin{displaymath}
\widehat{F}:=(I+V^{\ast}F_{1}F_{1}^{\ast}V+V^{\ast}F_{3}F_{3}^{\ast}V)^{-1}, \quad  \widehat{E}:=(I+U^{\ast}E_{1}^{\ast}E_{1}U+U^{\ast}E_{2}^{\ast}E_{2}U)^{-1}.
\end{displaymath}
We then have
\begin{displaymath}
(I+L^{\ast})(I+LL^{\ast})^{-1}Y^{-1}=\begin{pmatrix}
L_{1}&L_{3}\\
L_{2}&L_{4}\\
\end{pmatrix},\quad X^{-1}(I+R^{\ast}R)^{-1}(I+R^{\ast})=\begin{pmatrix}
R_{1}&R_{3}\\
R_{2}&R_{4}\\
\end{pmatrix},
\end{displaymath}
where 
\begin{align*}
L_{1}&=I-F_{1}^{\ast}V\widehat{F}V^{\ast},\quad L_{2}=-(I+F_{3}^{\ast}V)\widehat{F}V^{\ast},\quad L_{3}=F_{1}^{\ast}V\widehat{F},\quad L_{4}=(I+F_{3}^{\ast}V)\widehat{F},\\ 
R_{1}&=I-U\widehat{E}U^{\ast}E_{1}^{\ast},\quad R_{2}=-\widehat{E}U^{\ast}E_{1}^{\ast},\quad R_{3}=U\widehat{E}(I-U^{\ast}E_{2}^{\ast}),\quad R_{4}=\widehat{E}(I-U^{\ast}E_{2}^{\ast}).
\end{align*}
An application of Lemma 2.1 yields
\begin{displaymath}
M^{\dagger}=\begin{pmatrix}
(L_{1}N_{1}+L_{3}N_{2})R_{1}+(L_{1}N_{3}+L_{3}N_{4})R_{2}&(L_{1}N_{1}+L_{3}N_{2})R_{3}+(L_{1}N_{3}+L_{3}N_{4})R_{4}\\
(L_{2}N_{1}+L_{4}N_{2})R_{1}+(L_{2}N_{3}+L_{4}N_{4})R_{2}&(L_{2}N_{1}+L_{4}N_{2})R_{3}+(L_{2}N_{3}+L_{4}N_{4})R_{4}\\
\end{pmatrix}.
\end{displaymath}
Because
\begin{displaymath}
M=\begin{pmatrix}
I&-U\\
0&I\\
\end{pmatrix}\begin{pmatrix}
A&U\\
-V^{\ast}&I\\
\end{pmatrix}\begin{pmatrix}
I&0\\
V^{\ast}&I\\
\end{pmatrix}=\begin{pmatrix}
A+UV^{\ast}&0\\
0&I\\
\end{pmatrix},
\end{displaymath}
we immediately have
\begin{displaymath}
M^{\dagger}=\begin{pmatrix}
A+UV^{\ast}&0\\
0&I\\
\end{pmatrix}^{\dagger}=\begin{pmatrix}
(A+UV^{\ast})^{\dagger}&0\\
0&I\\
\end{pmatrix}.
\end{displaymath}
Due to $(A+UV^{\ast})^{\dagger}\in\mathbb{C}^{n\times m}$ and $(L_{1}N_{1}+L_{3}N_{2})R_{1}+(L_{1}N_{3}+L_{3}N_{4})R_{2}\in\mathbb{C}^{n\times m}$, it follows that
\begin{align*}
(A+UV^{\ast})^{\dagger}=(L_{1}N_{1}+L_{3}N_{2})R_{1}+(L_{1}N_{3}+L_{3}N_{4})R_{2}.\tag{3.19}
\end{align*}

(ii) Nevertheless, the expression (3.19) is not very legible. We next devote to simplifying (3.19). From $E_{N}N=0$ and $NF_{N}=0$, we have
\begin{displaymath}
E_{1}U+E_{3}=0,\quad E_{2}U+E_{4}=0, \quad -V^{\ast}F_{1}+F_{2}=0, \quad -V^{\ast}F_{3}+F_{4}=0.
\end{displaymath}
Since both $E_{N}$ and $F_{N}$ are orthogonal projectors, it follows that
\begin{align*}
 E_{1}^{\ast}=E_{1},\quad E_{3}^{\ast}&=E_{2},\quad E_{4}^{\ast}=E_{4},\quad
F_{1}^{\ast}=F_{1},\quad F_{2}^{\ast}=F_{3},\quad F_{4}^{\ast}=F_{4},\\
&E_{2}E_{3}+E_{4}^{2}=E_{4},\quad F_{2}F_{3}+F_{4}^{2}=F_{4}.
\end{align*}
Hence,
\begin{align*}
\widehat{F}&=(I+V^{\ast}F_{1}F_{1}^{\ast}V+V^{\ast}F_{3}F_{3}^{\ast}V)^{-1}=(I+F_{2}F_{2}^{\ast}+F_{4}F_{4}^{\ast})^{-1}\\
&=(I+F_{2}F_{3}+F_{4}^{2})^{-1}=(I+F_{4})^{-1},\\ \widehat{E}&=(I+U^{\ast}E_{1}^{\ast}E_{1}U+U^{\ast}E_{2}^{\ast}E_{2}U)^{-1}=(I+E_{3}^{\ast}E_{3}+E_{4}^{\ast}E_{4})^{-1}\\
&=(I+E_{2}E_{3}+E_{4}^{2})^{-1}=(I+E_{4})^{-1}.
\end{align*}
By (3.3a), (3.16), and (3.17), we have 
\begin{displaymath}
E_{4}=I+V^{\ast}N_{3}-N_{4}=E_{S_{A}}(I+HH^{\ast}E_{S_{A}})^{-1},
\end{displaymath}
which gives
\begin{align*}
\widehat{E}&=(I+HH^{\ast}E_{S_{A}})(I+E_{S_{A}}+HH^{\ast}E_{S_{A}})^{-1}\\
&=I-E_{S_{A}}(I+E_{S_{A}}+HH^{\ast}E_{S_{A}})^{-1}\\
&=I-(I+E_{S_{A}}+E_{S_{A}}HH^{\ast})^{-1}E_{S_{A}}.\tag{3.20}
\end{align*}
By (3.3c), (3.3d), and (3.4b), we have
\begin{displaymath}
F_{2}-\widetilde{F}_{2}=-(I-N_{4}-N_{2}U)V^{\ast}=-(I+F_{S_{A}}K^{\ast}K)^{-1}F_{S_{A}}V^{\ast},
\end{displaymath}
which yields 
\begin{displaymath}
F_{4}=I-N_{2}U-N_{4}=(I+F_{S_{A}}K^{\ast}K)^{-1}F_{S_{A}}.
\end{displaymath}
Thus,
\begin{align*}
\widehat{F}&=(I+F_{S_{A}}+F_{S_{A}}K^{\ast}K)^{-1}(I+F_{S_{A}}K^{\ast}K)\\
&=I-(I+F_{S_{A}}+F_{S_{A}}K^{\ast}K)^{-1}F_{S_{A}}\\
&=I-F_{S_{A}}(I+F_{S_{A}}+K^{\ast}KF_{S_{A}})^{-1}.\tag{3.21}
\end{align*}
Owing to $E_{S_{A}}U^{\ast}=0$ and $VF_{S_{A}}=0$, by (3.20) and (3.21), we have $\widehat{E}U^{\ast}=U^{\ast}$ and $V\widehat{F}=V$. Hence,
\begin{align*}
L_{1}&=I-F_{1}^{\ast}V\widehat{F}V^{\ast}=I-F_{1}VV^{\ast}=I-F_{2}^{\ast}V^{\ast}=I-F_{3}V^{\ast},\tag{3.22a}\\ 
L_{3}&=F_{1}^{\ast}V\widehat{F}=F_{1}V=F_{2}^{\ast}=F_{3},\tag{3.22b}\\ 
R_{1}&=I-U\widehat{E}U^{\ast}E_{1}^{\ast}=I-UU^{\ast}E_{1}=I+UE_{3}^{\ast}=I+UE_{2},\tag{3.22c}\\
R_{2}&=-\widehat{E}U^{\ast}E_{1}^{\ast}=-U^{\ast}E_{1}=E_{3}^{\ast}=E_{2},\tag{3.22d}
\end{align*}
due to $F_{1}^{\ast}=F_{1}$, $F_{2}^{\ast}=F_{3}$,
$E_{1}^{\ast}=E_{1}$,  $E_{3}^{\ast}=E_{2}$,  $V^{\ast}F_{1}=F_{2}$, and $E_{1}U=-E_{3}$. Using $N^{\dagger}E_{N}=0$ and $F_{N}N^{\dagger}=0$, we obtain
\begin{align*}
N_{1}E_{1}+N_{3}E_{2}=0, \quad F_{1}N_{1}+F_{3}N_{2}=0, \quad  F_{1}N_{3}+F_{3}N_{4}=0.\tag{3.23}
\end{align*}
By substituting (3.22a)--(3.22d) into (3.19), we obtain from (3.23) that
\begin{align*}
(A+UV^{\ast})^{\dagger}&=\big[(I-F_{3}V^{\ast})N_{1}+F_{3}N_{2}\big](I+UE_{2})+\big[(I-F_{3}V^{\ast})N_{3}+F_{3}N_{4}\big]E_{2}\\
&=\big[(I-F_{3}V^{\ast})N_{1}-F_{1}N_{1}\big](I+UE_{2})+\big[(I-F_{3}V^{\ast})N_{3}-F_{1}N_{3}\big]E_{2}\\
&=(I-F_{3}V^{\ast}-F_{1})N_{1}(I+UE_{2})+(I-F_{3}V^{\ast}-F_{1})N_{3}E_{2}\\
&=(I-F_{3}V^{\ast}-F_{1})\big[N_{1}(I+UE_{2})-N_{1}E_{1}\big]\\
&=(I-F_{3}V^{\ast}-F_{1})N_{1}(I+UE_{2}-E_{1})\\
&=(N_{1}UV^{\ast}+N_{1}A)N_{1}(UV^{\ast}N_{1}+AN_{1})\\
&=N_{1}(A+UV^{\ast})N_{1}(A+UV^{\ast})N_{1},\tag{3.24}
\end{align*}
where 
\begin{displaymath}
N_{1}=\big(I+A^{\dagger}UF_{S_{A}}U^{\ast}(A^{\dagger})^{\ast}\big)^{-1}\big(A^{\dagger}-A^{\dagger}US_{A}^{\dagger}V^{\ast}A^{\dagger}\big)\big(I+(A^{\dagger})^{\ast}VE_{S_{A}}V^{\ast}A^{\dagger}\big)^{-1}.
\end{displaymath} 

(iii) Subsequently, we further simplify the expression (3.24). Because $V^{\ast}F_{A}=0$ and $E_{S_{A}}S_{A}=0$, it follows that
\begin{align*}
\big(I+(A^{\dagger})^{\ast}VE_{S_{A}}V^{\ast}A^{\dagger}\big)(A+UV^{\ast})&=A+UV^{\ast}+(A^{\dagger})^{\ast}VE_{S_{A}}V^{\ast}A^{\dagger}A+(A^{\dagger})^{\ast}VE_{S_{A}}V^{\ast}A^{\dagger}UV^{\ast}\\
&=A+UV^{\ast}-(A^{\dagger})^{\ast}VE_{S_{A}}V^{\ast}F_{A}+(A^{\dagger})^{\ast}VE_{S_{A}}S_{A}V^{\ast}\\
&=A+UV^{\ast}.
\end{align*} 
Similarly, by $E_{A}U=0$ and $S_{A}F_{S_{A}}=0$, we can derive
\begin{displaymath}
(A+UV^{\ast})\big(I+A^{\dagger}UF_{S_{A}}U^{\ast}(A^{\dagger})^{\ast}\big)=A+UV^{\ast}.
\end{displaymath}
Therefore,
\begin{align*}
\big(I+(A^{\dagger})^{\ast}VE_{S_{A}}V^{\ast}A^{\dagger}\big)^{-1}(A+UV^{\ast})&=A+UV^{\ast},\tag{3.25}\\
(A+UV^{\ast})\big(I+A^{\dagger}UF_{S_{A}}U^{\ast}(A^{\dagger})^{\ast}\big)^{-1}&=A+UV^{\ast}.\tag{3.26}
\end{align*}
Furthermore, by $E_{A}U=0$ and $UE_{S_{A}}=0$, we have
\begin{align*}
(A+UV^{\ast})\big(A^{\dagger}-A^{\dagger}US_{A}^{\dagger}V^{\ast}A^{\dagger}\big)&=AA^{\dagger}+UV^{\ast}A^{\dagger}-AA^{\dagger}US_{A}^{\dagger}V^{\ast}A^{\dagger}-UV^{\ast}A^{\dagger}US_{A}^{\dagger}V^{\ast}A^{\dagger}\\
&=AA^{\dagger}+UV^{\ast}A^{\dagger}+E_{A}US_{A}^{\dagger}V^{\ast}A^{\dagger}-US_{A}S_{A}^{\dagger}V^{\ast}A^{\dagger}\\
&=AA^{\dagger}+UE_{S_{A}}V^{\ast}A^{\dagger}+E_{A}US_{A}^{\dagger}V^{\ast}A^{\dagger}\\
&=AA^{\dagger}.\tag{3.27}
\end{align*}
Analogously, using $V^{\ast}F_{A}=0$ and 
$F_{S_{A}}V^{\ast}=0$, we obtain
\begin{align*}
\big(A^{\dagger}-A^{\dagger}US_{A}^{\dagger}V^{\ast}A^{\dagger}\big)(A+UV^{\ast})=A^{\dagger}A.\tag{3.28}
\end{align*}

Consequently, using (3.24)--(3.28), we arrive at
\begin{displaymath}
(A+UV^{\ast})^{\dagger}=\big(I+A^{\dagger}UF_{S_{A}}U^{\ast}(A^{\dagger})^{\ast}\big)^{-1}\big(A^{\dagger}-A^{\dagger}US_{A}^{\dagger}V^{\ast}A^{\dagger}\big)\big(I+(A^{\dagger})^{\ast}VE_{S_{A}}V^{\ast}A^{\dagger}\big)^{-1}.
\end{displaymath}
This completes the proof. \qed

\medskip

\noindent{\bf Remark 3.3.} Due to $A^{\dagger}UF_{S_{A}}U^{\ast}(A^{\dagger})^{\ast}$ and $(A^{\dagger})^{\ast}VE_{S_{A}}V^{\ast}A^{\dagger}$ are (Hermitian) positive semidefinite, it follows that $I+A^{\dagger}UF_{S_{A}}U^{\ast}(A^{\dagger})^{\ast}$ and $I+(A^{\dagger})^{\ast}VE_{S_{A}}V^{\ast}A^{\dagger}$ are (Hermitian) positive definite and hence they are nonsingular. Moreover, the factors $\left(I+A^{\dagger}UF_{S_{A}}U^{\ast}(A^{\dagger})^{\ast}\right)^{-1}$ and $\left(I+(A^{\dagger})^{\ast}VE_{S_{A}}V^{\ast}A^{\dagger}\right)^{-1}$ can be computed via the SMW formula (1.1), that is,
\begin{align*}
\big(I+A^{\dagger}UF_{S_{A}}U^{\ast}(A^{\dagger})^{\ast}\big)^{-1}&=I-A^{\dagger}U\big(I+F_{S_{A}}U^{\ast}(A^{\dagger})^{\ast}A^{\dagger}U\big)^{-1}F_{S_{A}}U^{\ast}(A^{\dagger})^{\ast},\\
\big(I+(A^{\dagger})^{\ast}VE_{S_{A}}V^{\ast}A^{\dagger}\big)^{-1}&=I-(A^{\dagger})^{\ast}V\big(I+E_{S_{A}}V^{\ast}A^{\dagger}(A^{\dagger})^{\ast}V\big)^{-1}E_{S_{A}}V^{\ast}A^{\dagger}.
\end{align*}
For a given matrix $A\in\mathbb{C}^{m\times n}$, we assume that $A^{\dagger}$ has been precomputed. For a varied perturbation $UV^{\ast}$ (here $U\in\mathbb{C}^{m\times r}$ and $V\in\mathbb{C}^{n\times r}$), we can compute $(A+UV^{\ast})^{\dagger}$ via the singular value decomposition of $A+UV^{\ast}$. However, this approach is expensive, especially when $m$ and $n$ are large. Fortunately, if $r$ is much smaller than $m$ and $n$, $S_{A}^{\dagger}=\left(I+V^{\ast}A^{\dagger}U\right)^{\dagger}$ is much easier to compute than $(A+UV^{\ast})^{\dagger}$. Therefore, the formula (3.18) provides an effective method to compute $(A+UV^{\ast})^{\dagger}$, provided that the conditions of Theorem 3.2 are satisfied.

\medskip

On the basis of Theorem 3.2, we can derive some simple expressions for $(A+UV^{\ast})^{\dagger}$ under several special conditions.

\medskip

\noindent{\bf Corollary 3.4.} \emph{If $\mathcal{R}(U)\subseteq \mathcal{R}(A)$, $\mathcal{R}(V)\subseteq \mathcal{R}(A^{\ast})$, $\mathcal{R}(U^{\ast})\subseteq \mathcal{R}(S_{A})\cap \mathcal{R}(S_{A}^{\ast})$, and $\mathcal{R}(V^{\ast})\subseteq \mathcal{R}(S_{A})\cap \mathcal{R}(S_{A}^{\ast})$, then} 
\begin{displaymath}	
(A+UV^{\ast})^{\dagger}=A^{\dagger}-A^{\dagger}U(I+V^{\ast}A^{\dagger}U)^{\dagger}V^{\ast}A^{\dagger}.
\end{displaymath}

\noindent{\bf Corollary 3.5.} \emph{If $\mathcal{R}(U)\subseteq \mathcal{R}(A)$, $\mathcal{R}(V)\subseteq \mathcal{R}(A^{\ast})$, and $S_{A}=I+V^{\ast}A^{\dagger}U$ is nonsingular, then} 
\begin{displaymath}
(A+UV^{\ast})^{\dagger}=A^{\dagger}-A^{\dagger}U(I+V^{\ast}A^{\dagger}U)^{-1}V^{\ast}A^{\dagger}.
\end{displaymath}

\noindent{\bf Corollary 3.6.} \emph{Let $A\in\mathbb{C}^{n\times n}$, $U\in\mathbb{C}^{n\times r}$, and $V\in\mathbb{C}^{n\times r}$. If $A$ and $S_{A}=I+V^{\ast}A^{-1}U$ are nonsingular, then 
\begin{displaymath}
(A+UV^{\ast})^{\dagger}=A^{-1}-A^{-1}U(I+V^{\ast}A^{-1}U)^{-1}V^{\ast}A^{-1}.
\end{displaymath}
Indeed, $A+UV^{\ast}$ is also nonsingular and $(A+UV^{\ast})^{-1}=A^{-1}-A^{-1}U(I+V^{\ast}A^{-1}U)^{-1}V^{\ast}A^{-1}$, which is exactly the SMW formula {\rm (1.1)}.
}

\medskip

Under the conditions $\mathcal{R}(U)\subseteq \mathcal{R}(A)$ and $\mathcal{R}(V)\subseteq \mathcal{R}(A^{\ast})$, we next give a necessary and sufficient condition to validate (3.18).

\medskip

\noindent{\bf Theorem 3.7.} \emph{Let $A\in \mathbb{C}^{m\times n}$, $U\in \mathbb{C}^{m\times r}$, and $V\in \mathbb{C}^{n\times r}$. If $\mathcal{R}(U)\subseteq \mathcal{R}(A)$ and $\mathcal{R}(V)\subseteq \mathcal{R}(A^{\ast})$, then {\rm (3.18)} holds if and only if both $UE_{S_{A}}V^{\ast}A^{\dagger}$ and $A^{\dagger}UF_{S_{A}}V^{\ast}$ are Hermitian and $A^{\dagger}UF_{S_{A}}E_{S_{A}}V^{\ast}A^{\dagger}=0$, where $S_{A}=I+V^{\ast}A^{\dagger}U$.}

\medskip

\noindent{\bf Proof.} According to the derivations of (3.25) and (3.26), by $\mathcal{R}(U)\subseteq \mathcal{R}(A)$ and $\mathcal{R}(V)\subseteq \mathcal{R}(A^{\ast})$, we obtain that (3.25) and (3.26) hold.

(i) Let $\widehat{A}=\big(I+A^{\dagger}UF_{S_{A}}U^{\ast}(A^{\dagger})^{\ast}\big)^{-1}\big(A^{\dagger}-A^{\dagger}US_{A}^{\dagger}V^{\ast}A^{\dagger}\big)\big(I+(A^{\dagger})^{\ast}VE_{S_{A}}V^{\ast}A^{\dagger}\big)^{-1}$.
Because $E_{A}U=0$, $V^{\ast}F_{A}=0$, and $E_{S_{A}}S_{A}=0$, it follows that
\begin{align*}
(A+UV^{\ast})\widehat{A}(A+UV^{\ast})&=(A+UV^{\ast})\big(A^{\dagger}-A^{\dagger}US_{A}^{\dagger}V^{\ast}A^{\dagger}\big)(A+UV^{\ast})\\
&=\big(AA^{\dagger}+UE_{S_{A}}V^{\ast}A^{\dagger}\big)(A+UV^{\ast})\\
&=A+AA^{\dagger}UV^{\ast}+UE_{S_{A}}V^{\ast}A^{\dagger}A+UE_{S_{A}}V^{\ast}A^{\dagger}UV^{\ast}\\
&=A+UV^{\ast}+UE_{S_{A}}V^{\ast}+UE_{S_{A}}(S_{A}-I)V^{\ast}\\
&=A+UV^{\ast},\tag{3.29}
\end{align*}
where we have used the equalities (3.25) and (3.26). 

(ii) By (3.25) and (3.26), we can easily get that
\begin{displaymath}
\widehat{A}(A+UV^{\ast})\widehat{A}=\big(I+A^{\dagger}UF_{S_{A}}U^{\ast}(A^{\dagger})^{\ast}\big)^{-1}\widetilde{A}\big(I+(A^{\dagger})^{\ast}VE_{S_{A}}V^{\ast}A^{\dagger}\big)^{-1},
\end{displaymath}
where 
\begin{displaymath}
\widetilde{A}=\big(A^{\dagger}-A^{\dagger}US_{A}^{\dagger}V^{\ast}A^{\dagger}\big)(A+UV^{\ast})\big(A^{\dagger}-A^{\dagger}US_{A}^{\dagger}V^{\ast}A^{\dagger}\big).
\end{displaymath}
Due to $E_{A}U=0$ and $S_{A}^{\dagger}E_{S_{A}}=0$, it follows that
\begin{align*}
\widetilde{A}&=\big(A^{\dagger}-A^{\dagger}US_{A}^{\dagger}V^{\ast}A^{\dagger}\big)\big(AA^{\dagger}+UE_{S_{A}}V^{\ast}A^{\dagger}\big)\\
&=A^{\dagger}-A^{\dagger}US_{A}^{\dagger}V^{\ast}A^{\dagger}+\big(A^{\dagger}-A^{\dagger}US_{A}^{\dagger}V^{\ast}A^{\dagger}\big)UE_{S_{A}}V^{\ast}A^{\dagger}\\
&=A^{\dagger}-A^{\dagger}US_{A}^{\dagger}V^{\ast}A^{\dagger}+\big(A^{\dagger}UE_{S_{A}}V^{\ast}A^{\dagger}-A^{\dagger}US_{A}^{\dagger}(S_{A}-I)E_{S_{A}}V^{\ast}A^{\dagger}\big)\\
&=A^{\dagger}-A^{\dagger}US_{A}^{\dagger}V^{\ast}A^{\dagger}+A^{\dagger}UF_{S_{A}}E_{S_{A}}V^{\ast}A^{\dagger}.
\end{align*}
Hence, $\widehat{A}(A+UV^{\ast})\widehat{A}=\widehat{A}$ if and only if $A^{\dagger}UF_{S_{A}}E_{S_{A}}V^{\ast}A^{\dagger}=0$. 

(iii) Using (3.25), (3.26), $E_{A}U=0$, and $V^{\ast}F_{A}=0$, we obtain
\begin{align*}
(A+UV^{\ast})\widehat{A}&=\big(AA^{\dagger}+UE_{S_{A}}V^{\ast}A^{\dagger}\big)\big(I+(A^{\dagger})^{\ast}VE_{S_{A}}V^{\ast}A^{\dagger}\big)^{-1},\\
\widehat{A}(A+UV^{\ast})&=\big(I+A^{\dagger}UF_{S_{A}}U^{\ast}(A^{\dagger})^{\ast}\big)^{-1}\big(A^{\dagger}A+A^{\dagger}UF_{S_{A}}V^{\ast}\big).
\end{align*}
Then we have
\begin{align*}
\big((A+UV^{\ast})\widehat{A}\big)^{\ast}&=\big(I+(A^{\dagger})^{\ast}VE_{S_{A}}V^{\ast}A^{\dagger}\big)^{-1}\big(AA^{\dagger}+(UE_{S_{A}}V^{\ast}A^{\dagger})^{\ast}\big),\\
\big(\widehat{A}(A+UV^{\ast})\big)^{\ast}&=\big(A^{\dagger}A+(A^{\dagger}UF_{S_{A}}V^{\ast})^{\ast}\big)\big(I+A^{\dagger}UF_{S_{A}}U^{\ast}(A^{\dagger})^{\ast}\big)^{-1}.
\end{align*}
Note that 
\begin{align*}
\big(I+(A^{\dagger})^{\ast}VE_{S_{A}}V^{\ast}A^{\dagger}\big)^{-1}AA^{\dagger}&=AA^{\dagger}\big(I+(A^{\dagger})^{\ast}VE_{S_{A}}V^{\ast}A^{\dagger}\big)^{-1},\\
\big(I+A^{\dagger}UF_{S_{A}}U^{\ast}(A^{\dagger})^{\ast}\big)^{-1}A^{\dagger}A&=A^{\dagger}A\big(I+A^{\dagger}UF_{S_{A}}U^{\ast}(A^{\dagger})^{\ast}\big)^{-1}.
\end{align*}
Hence, $(A+UV^{\ast})\widehat{A}$ is Hermitian if and only if 
$UE_{S_{A}}V^{\ast}A^{\dagger}\left(I+(A^{\dagger})^{\ast}VE_{S_{A}}V^{\ast}A^{\dagger}\right)^{-1}$ is Hermitian, and $\widehat{A}(A+UV^{\ast})$ is Hermitian if and only if $\left(I+A^{\dagger}UF_{S_{A}}U^{\ast}(A^{\dagger})^{\ast}\right)^{-1} A^{\dagger}UF_{S_{A}}V^{\ast}$ is Hermitian, that is,
\begin{align*}
\big(I+(A^{\dagger})^{\ast}VE_{S_{A}}V^{\ast}A^{\dagger}\big)UE_{S_{A}}V^{\ast}A^{\dagger}&=(UE_{S_{A}}V^{\ast}A^{\dagger})^{\ast}\big(I+(A^{\dagger})^{\ast}VE_{S_{A}}V^{\ast}A^{\dagger}\big),\tag{3.30}\\
A^{\dagger}UF_{S_{A}}V^{\ast}\big(I+A^{\dagger}UF_{S_{A}}U^{\ast}(A^{\dagger})^{\ast}\big)&=\big(I+A^{\dagger}UF_{S_{A}}U^{\ast}(A^{\dagger})^{\ast}\big)(A^{\dagger}UF_{S_{A}}V^{\ast})^{\ast}.\tag{3.31}
\end{align*}
From (3.30) and (3.31), we obtain that  $(A+UV^{\ast})\widehat{A}$ is Hermitian if and only if $UE_{S_{A}}V^{\ast}A^{\dagger}$ is Hermitian, and $\widehat{A}(A+UV^{\ast})$ is Hermitian if and only if $A^{\dagger}UF_{S_{A}}V^{\ast}$ is Hermitian.

According to the definition of the MP inverse of a matrix, it follows that the statement in Theorem 3.7 holds. \qed

\medskip

\noindent{\bf Remark 3.8.} Under the conditions $\mathcal{R}(U)\subseteq \mathcal{R}(A)$ and $\mathcal{R}(V)\subseteq \mathcal{R}(A^{\ast})$, by comparing Theorem 3.2 with Theorem 3.7, we can readily observe that $\mathcal{R}(U^{\ast})\subseteq \mathcal{R}(S_{A})$ and $\mathcal{R}(V^{\ast})\subseteq \mathcal{R}(S_{A}^{\ast})$ are only sufficient conditions to validate (3.18). In addition, from the proof of (3.29), we find that
\begin{displaymath}
\big(I+A^{\dagger}UF_{S_{A}}U^{\ast}(A^{\dagger})^{\ast}\big)^{-1}\big(A^{\dagger}-A^{\dagger}US_{A}^{\dagger}V^{\ast}A^{\dagger}\big)\big(I+(A^{\dagger})^{\ast}VE_{S_{A}}V^{\ast}A^{\dagger}\big)^{-1}
\end{displaymath}
is a $\{1\}$-inverse of $A+UV^{\ast}$, provided that $\mathcal{R}(U)\subseteq \mathcal{R}(A)$ and $\mathcal{R}(V)\subseteq \mathcal{R}(A^{\ast})$.

\medskip
\bigskip

\noindent{\bf \large 4. Extensions}

\medskip

In this section, we devote to extending the results stated in Theorems 3.2 and 3.7 to the bounded linear operators case.

Let $\mathcal{H}_{1}$ and $\mathcal{H}_{2}$ be two Hilbert spaces over the same field. The set of all bounded linear operators from $\mathcal{H}_{1}$ into $\mathcal{H}_{2}$ is denoted by $\mathcal{B}(\mathcal{H}_{1},\mathcal{H}_{2})$. For any $T\in\mathcal{B}(\mathcal{H}_{1},\mathcal{H}_{2})$, let $T^{\ast}$, $\mathcal{R}(T)$, and $\mathcal{N}(T)$ denote the adjoint, the range, and the null space of $T$, respectively. The Moore--Penrose inverse of $T\in\mathcal{B}(\mathcal{H}_{1},\mathcal{H}_{2})$ is denoted by $T^{\dagger}\in\mathcal{B}(\mathcal{H}_{2},\mathcal{H}_{1})$ (if it exists), which is defined as the unique operator $Z\in\mathcal{B}(\mathcal{H}_{2},\mathcal{H}_{1})$ satisfying the four equations
\begin{align*}
TZT=T, \quad ZTZ=Z, \quad  (TZ)^{\ast}=TZ, \quad (ZT)^{\ast}=ZT.\tag{4.1}
\end{align*}
Unlike the matrices case, $T^{\dagger}$ is not always existent. Indeed, an operator $T$ has the MP inverse if and only if $\mathcal{R}(T)$ is closed [17].

In view of Theorems 3.2 and 3.7, we can establish the following more general results whose detailed proofs are omitted, because we can directly check that the operator 
\begin{displaymath}
\big(I+A^{\dagger}UF_{S_{A}}U^{\ast}(A^{\dagger})^{\ast}\big)^{-1}\big(A^{\dagger}-A^{\dagger}US_{A}^{\dagger}V^{\ast}A^{\dagger}\big)\big(I+(A^{\dagger})^{\ast}VE_{S_{A}}V^{\ast}A^{\dagger}\big)^{-1}
\end{displaymath}
satisfies the four equations in (4.1).

\medskip

\noindent{\bf Theorem 4.1.} \emph{Let $A\in\mathcal{B}(\mathcal{H}_{1},\mathcal{H}_{2})$, $U\in\mathcal{B}(\mathcal{H}_{3},\mathcal{H}_{2})$, and $V\in\mathcal{B}(\mathcal{H}_{3},\mathcal{H}_{1})$, where $\mathcal{H}_{i}$ $(i=1,2,3)$ are Hilbert spaces over the same field. Assume that $\mathcal{R}(A)$ and $\mathcal{R}(S_{A})$ are closed. If $\mathcal{R}(U)\subseteq \mathcal{R}(A)$, $\mathcal{R}(V)\subseteq \mathcal{R}(A^{\ast})$, $\mathcal{R}(U^{\ast})\subseteq \mathcal{R}(S_{A})$, and $\mathcal{R}(V^{\ast})\subseteq \mathcal{R}(S_{A}^{\ast})$, then $\mathcal{R}(A+UV^{\ast})$ is closed and 
\begin{align*}
(A+UV^{\ast})^{\dagger}=\big(I+A^{\dagger}UF_{S_{A}}U^{\ast}(A^{\dagger})^{\ast}\big)^{-1}\big(A^{\dagger}-A^{\dagger}US_{A}^{\dagger}V^{\ast}A^{\dagger}\big)\big(I+(A^{\dagger})^{\ast}VE_{S_{A}}V^{\ast}A^{\dagger}\big)^{-1},\tag{4.2}
\end{align*}
where $S_{A}=I+V^{\ast}A^{\dagger}U$, $E_{S_{A}}=I-S_{A}S_{A}^{\dagger}$, and $F_{S_{A}}=I-S_{A}^{\dagger}S_{A}$.} 

\bigskip

\noindent{\bf Theorem 4.2.} \emph{Let $A\in\mathcal{B}(\mathcal{H}_{1},\mathcal{H}_{2})$, $U\in\mathcal{B}(\mathcal{H}_{3},\mathcal{H}_{2})$, and $V\in\mathcal{B}(\mathcal{H}_{3},\mathcal{H}_{1})$, where $\mathcal{H}_{i}$ $(i=1,2,3)$ are Hilbert spaces over the same field. Assume that $\mathcal{R}(A)$, $\mathcal{R}(S_{A})$, and $\mathcal{R}(A+UV^{\ast})$ are closed. If $\mathcal{R}(U)\subseteq \mathcal{R}(A)$, $\mathcal{R}(V)\subseteq \mathcal{R}(A^{\ast})$, then {\rm (4.2)} holds if and only if both $UE_{S_{A}}V^{\ast}A^{\dagger}$ and $A^{\dagger}UF_{S_{A}}V^{\ast}$ are self-adjoint and $A^{\dagger}UF_{S_{A}}E_{S_{A}}V^{\ast}A^{\dagger}=0$, where $S_{A}=I+V^{\ast}A^{\dagger}U$.
}

\medskip
\bigskip

\noindent{\bf \large Acknowledgements}

\medskip

The author would like to thank the anonymous referee for his/her valuable comments and suggestions, which greatly improved the original manuscript of this paper. The author is grateful to Professor Chen-Song Zhang for his helpful suggestions. This work was supported partially by the Science Challenge Project (Grant No. JCKY201612A503-1-303) and the National Natural Science Foundation of China (Grant Nos. 91430215, 91530323).

\medskip
\bigskip

\noindent{\bf \large References}

\medskip

\small
{
	
[1] A. Ben-Israel, T.N.E. Greville, Generalized Inverses: Theory and Applications, 2nd edition, Springer-Verlag, New York, 2003.

[2] J. Sherman, W.J. Morrison, Adjustment of an inverse matrix corresponding to a change in one element of a given matrix, Ann. Math. Stat. 21 (1950)
124--127.

[3] M.A. Woodbury, Inverting modified matrices, Memorandum Report 42, Statistical Research Group, Princeton University, Princeton, NJ, 1950.

[4] M.S. Bartlett, An inverse matrix adjustment arising in discriminant analysis, Ann. Math. Stat. 22 (1951) 107--111.

[5] W.W. Hager, Updating the inverse of a matrix, SIAM Rev. 31 (1989) 221--239.

[6] A. Malyshev, M. Sadkane, Using the Sherman–Morrison–Woodbury inversion formula for a fast solution of tridiagonal block Toeplitz systems, Linear Algebra Appl. 435 (2011) 2693--2707.

[7] S.-H. Lai, B.C. Vemuri, Sherman--Morrison--Woodbury--Formula--Based Algorithms for the Surface Smoothing Problem, Linear Algebra Appl. 265 (1997) 203--229.

[8] \AA{}. Bj\"{o}rck, Numerical Methods for Least Squares Problems, Society for Industrial and Applied
Mathematics, Philadelphia, PA, 1996.

[9] G.W. Stewart, J.-G. Sun, Matrix Perturbation Theory, Academic Press, New York, 1990.

[10] P.-\AA{}. Wedin, Perturbation theory for pseudo-inverses, BIT 13 (1973) 217--232.

[11] N. Castro-Gonz\'{a}lez, F.M. Dopico, J.M. Molera, Multiplicative perturbation theory of
the Moore--Penrose inverse and the least squares problem, Linear Algebra Appl. 503 (2016) 1--25.

[12] N. Castro-Gonz\'{a}lez, M.F. Mart\'{i}nez-Serrano, J. Robles, Expressions for the Moore--Penrose inverse of block matrices involving the Schur complement, Linear Algebra Appl. 471 (2015) 353--368.

[13] J.K. Baksalary, O.M. Baksalary, Particular formulae for the Moore--Penrose inverse of a columnwise
partitioned matrix, Linear Algebra Appl. 421 (2007) 16--23.

[14] D.S. Cvetkovi\'{c}-Ili\'{c}, Expression of the Drazin and MP-inverse of partitioned matrix and quotient
identity of generalized Schur complement, Appl. Math. Comput. 213 (2009) 18--24.

[15] C.Y. Deng, H.K. Du, Representations of the Moore--Penrose inverse for a class of 2-by-2 block
operator valued partial matrices, Linear Multilinear Algebra 58 (2010) 15--26.

[16] X. Liu, H. Jin, D.S. Cvetkovi\'{c}-Ili\'{c}, Representations of generalized inverses of partitioned matrix
involving Schur complement, Appl. Math. Comput. 219 (2013) 9615--9629.

[17] W.V. Petryshyn, On generalized inverses and on the uniform convergence of $(I-\beta K)^{n}$ with application to iterative methods, J. Math. Anal. Appl. 18 (1967), 417--439.

}

\end{document}